\documentclass[letterpaper, 10 pt, conference]{ieeeconf}  

\IEEEoverridecommandlockouts                              
\usepackage{mathtools}
\usepackage{amsmath,amssymb,amsfonts}
\usepackage{hyperref}
\usepackage{cancel}
\usepackage{array}
\usepackage{subcaption}
\DeclareCaptionLabelSeparator{periodspace}{.\quad}
\usepackage[nameinlink,capitalise]{cleveref}
\crefname{equation}{}{} 

\usepackage{standalone}
\usepackage{breakcites}
\usepackage{cite}
\usepackage{multicol}

\usepackage{tikz}
\usetikzlibrary{shapes,arrows,positioning,calc}
\usepackage{stanli}
\usetikzlibrary{decorations.pathreplacing}

\usepackage{pgfplots}
\pgfplotsset{compat=1.18}
\usepgfplotslibrary{colorbrewer} 
\usepgfplotslibrary{fillbetween}

\usepackage{pgfplotstable}
\usepackage{booktabs}
\usepackage{colortbl}

\newtheorem{definition}{Definition}

\newtheorem{theorem}{Theorem}

\newcommand*\circled[1]{\tikz[baseline=(char.base)]{\node(char)[shape=rounded rectangle,draw,inner sep=0.6pt,minimum height=1.5ex]{#1};}} %
\newcommand\kronF[2]{{#1}^{\circled{\tiny{\ensuremath{#2}}}}} 
\renewcommand\Vec[1][]{\textnormal{vec}\ensuremath{\if$#1$ \else \left[#1\right]\fi}}

\newcommand{\real}{\mathbb{R}}

\newcommand{\bzero}{\ensuremath{\mathbf{0}}} 

\newcommand{\tv}{\ensuremath{\tilde{v}}}
\newcommand{\tw}{\ensuremath{\tilde{w}}}

\newcommand{\bA}{\ensuremath{\mathbf{A}}}
\newcommand{\bB}{\ensuremath{\mathbf{B}}}
\newcommand{\bC}{\ensuremath{\mathbf{C}}}
\newcommand{\bD}{\ensuremath{\mathbf{D}}}
\newcommand{\bI}{\ensuremath{\mathbf{I}}}
\newcommand{\bM}{\ensuremath{\mathbf{M}}}
\newcommand{\bS}{\ensuremath{\mathbf{S}}}
\newcommand{\bV}{\ensuremath{\mathbf{V}}}
\newcommand{\bW}{\ensuremath{\mathbf{W}}}

\newcommand{\ba}{\ensuremath{\mathbf{a}}}

\renewcommand{\bf}{\ensuremath{\mathbf{f}}}
\newcommand{\bg}{\ensuremath{\mathbf{g}}}
\newcommand{\bh}{\ensuremath{\mathbf{h}}}
\newcommand{\bu}{\ensuremath{\mathbf{u}}}
\newcommand{\bv}{\ensuremath{\mathbf{v}}}
\newcommand{\bw}{\ensuremath{\mathbf{w}}}
\newcommand{\bx}{\ensuremath{\mathbf{x}}}
\newcommand{\by}{\ensuremath{\mathbf{y}}}

\newcommand{\cE}{\ensuremath{\mathcal{E}}}
\newcommand{\cH}{\ensuremath{\mathcal{H}}}
\newcommand{\cL}{\ensuremath{\mathcal{L}}}

\newcommand{\perm}[2]{\ensuremath{\bS_{#1 \times #2}}}

\newcommand{\F}{\ensuremath{\mathbf{F}}}

\title{\LARGE \textbf{
Scalable Computation of $\cH_\infty$ Energy Functions \\for Polynomial Drift Nonlinear Systems*
}}

\author{Nicholas A. Corbin$^{1}$ and Boris Kramer$^{2}$
\thanks{*This work was supported by the National Science Foundation under Grant CMMI-2130727.}
\thanks{$^{1}$N. Corbin is with the Department of Mechanical and Aerospace Engineering, University of California San Diego, La Jolla, CA 92093-0411 USA {\tt\small ncorbin@ucsd.edu}}
\thanks{$^{2}$B. Kramer is with the Faculty of the Department of Mechanical and Aerospace Engineering, University of California San Diego, La Jolla, CA 92093-0411 USA {\tt\small bmkramer@ucsd.edu}}
}

\begin{document}

\maketitle
\thispagestyle{empty}
\pagestyle{empty}

\begin{abstract}
        This paper presents a scalable tensor-based approach to computing controllability and observability-type energy functions for nonlinear dynamical systems with polynomial drift and linear input and output maps.
        Using Kronecker product polynomial expansions, we convert the Hamilton-Jacobi-Bellman partial differential equations for the energy functions into a series of
        algebraic equations for the coefficients of the energy functions.
        We derive the specific tensor structure that arises from the Kronecker product representation and analyze the computational complexity to efficiently solve these equations.
        The convergence and scalability of the proposed energy function computation approach is demonstrated on a nonlinear reaction-diffusion model with cubic drift nonlinearity, for which we compute degree 3 energy function approximations in $n=1023$ dimensions and degree 4 energy function approximations in $n=127$ dimensions.
\end{abstract}


\section{Introduction}\label{sec:intro}
In many modern engineering applications, models are both high-dimensional and nonlinear.
Common approaches to controller design often require sacrificing either the high-dimensionality or the nonlinearity in order to remain
computationally tractable.
For applications where both the nonlinearity and the high-dimensionality are critical to the model's accuracy, nonlinear model reduction is an attractive way to systematically reduce the system dimensionality---to keep computations manageable---while retaining important nonlinear features of the model.

Balanced truncation (BT) is one of the predominant system-theoretic model reduction approaches for linear time-invariant (LTI) systems \cite{Mullis1976,Moore1981,Gugercin2004}.
Its success for LTI systems has stimulated much interest in extending BT to nonlinear systems.
While Scherpen provided the theoretical extensions to nonlinear control-affine systems \cite{Scherpen1993,Scherpen1994,Scherpen1994a,Scherpen1996}, developing \emph{scalable computational methods} to implement nonlinear BT for nonlinear systems remains an active area of research.

One of the main computational challenges in nonlinear BT is solving the Hamilton-Jacobi-Bellman (HJB) partial differential equations (PDEs) for the controllability and observability-type nonlinear \emph{energy functions}.
HJB PDEs are notoriously difficult to solve for general nonlinear systems, \emph{especially} for high-dimensional models of interest in model reduction.
Solving HJB PDEs is also required in nonlinear optimal control problems, so many approaches have been developed to approximate solutions to HJB PDEs, including
state-dependent Riccati equations~\cite{Cimen2012}, algebraic Gramians \cite{Condon2005,Gray2006,Benner2017}, discretization techniques \cite{Falcone2016}, and iterative approaches \cite{Kalise2018,Dolgov2021}.
We are particularly interested in Al'brekht's method \cite{Albrekht1961,Lukes1969}, a power-series-based approach which has been widely used to locally approximate HJB PDE solutions for low-dimensional models \cite{Garrard1977,Scherpen1994a,Fujimoto2008a,Fujimoto2010,Krener2013,Sahyoun2013,Krener2019}.
More recently, this power series approach has been adapted to solve optimal control problems for moderately-sized
bilinear systems \cite{Breiten2018,Breiten2019}, quadratic drift systems \cite{Breiten2019a,Borggaard2020}, and polynomial drift systems \cite{Borggaard2021,Almubarak2019}.
These works demonstrated that polynomial approximation induces tensor structure in the resulting equations for the polynomial coefficients that modern solvers (e.g., \cite{Grasedyck2004,Chen2019}) can exploit to improve scalability.
For computing nonlinear BT energy functions, this approach has been adapted in \cite{Kramer2022a},
but only for models with quadratic drift nonlinearities.

The first contribution of this article is a scalable Kronecker product-based approach to computing nonlinear BT energy function approximations for systems with \emph{polynomial drift of arbitrary degree} $\ell$, linear inputs, and linear outputs.
Moreover, we provide the explicit form of the equations of the energy function coefficients in Kronecker product form.
Additionally, we provide open-access software implementations for the proposed algorithms in the \texttt{cnick1/NLBalancing} repository \cite{NLBalancing2023} under the \texttt{v0.9.1} tag, along with the code to reproduce the results for the presented example.

The rest of the article is organized as follows.
\Cref{sec:background} reviews the definitions for the $\cH_\infty$ energy functions used in nonlinear balancing, along with a summary of Al'brekht's method for approximating energy function solutions and an introduction to the Kronecker product notation.
\Cref{sec:NLBT-Poly} presents the main contribution of the article, namely the algorithms for computing the polynomial energy function approximations for polynomial drift systems.
In \cref{sec:complexity}, we discuss crucial details to efficiently implement the proposed method in a scalable manner.
We then demonstrate the method on a nonlinear heat equation finite element model in \cref{sec:fem}, before concluding the article in \cref{sec:conclusion}.


\section{Background, Definitions, and Notation}\label{sec:background}
We first review the definitions for the $\cH_\infty$ nonlinear~BT energy functions in \cref{sec:NLBT-background}.
In \cref{sec:albrekht}, we review Al'brekht's method for locally approximating solutions to HJB PDEs.
Afterwards, since our method is based on the Kronecker product polynomial representation, we review basic notation and definitions relating to Kronecker product polynomial expansions in \cref{sec:notation}.

\subsection{Energy Functions for \texorpdfstring{$\cH_\infty$}{H-infinity} Nonlinear Balancing}\label{sec:NLBT-background}
Consider the nonlinear control-affine dynamical system
\begin{equation}\label{eq:FOM-NL}
        \dot{\bx}(t)  = \bf(\bx(t))  + \bg(\bx(t)) \bu(t), \quad
        \by(t)        = \bh(\bx(t)),
\end{equation}
where $\bx(t) \in \real^n$ is the state, $\bf \colon \real^n \to \real^n$ is the nonlinear drift, $\bg \colon \real^n \to \real^{n \times m}$ is the nonlinear input map, $\bh \colon \real^n \to \real^{p}$ is the nonlinear output map, $\bu(t) \in \real^m$ is a vector of input signals, and $\by(t) \in \real^p$ is a vector of output signals.

The $\cH_\infty$ nonlinear balancing framework \cite{Scherpen1996} defines a pair of energy functions that generalize the concepts of controllability and observability to (potentially unstable) systems of the form \cref{eq:FOM-NL}.
These energy functions (defined next) are then balanced using a nonlinear state-space transformation, and model reduction involves truncating states determined to be less important in the balanced representation.
\begin{definition}\cite[Def. 5.1]{Scherpen1996}
        Let $\gamma$ be a positive constant $\gamma > 0, \gamma \neq 1$, and define $\eta \coloneqq 1-\gamma^{-2}$.
        The $\cH_\infty$ past energy of the nonlinear system \cref{eq:FOM-NL} is defined as
        \begin{equation} \label{eq:HinftyPastEnergyDef}
                \cE_\gamma^{-}(\bx_0)  \coloneqq \!\!\!\! \min_{\substack{\bu \in L_{2}(-\infty, 0] \\ \bx(-\infty) = \bzero,\,  \bx(0) = \bx_0}} \! \frac{1}{2} \int\displaylimits_{-\infty}^{0} \eta \Vert \by(t) \Vert^2  +  \Vert \bu(t) \Vert^2 {\rm{d}}t.
        \end{equation}
        If $\gamma<1$, the $\cH_\infty$ future energy of the nonlinear system \cref{eq:FOM-NL} is defined as
        \begin{equation} \label{eq:HinftyFutureEnergyDef-a}
                \cE_\gamma^{+}(\bx_0)  \coloneqq \!\!\!\! \max_{\substack{\bu \in L_{2}[0,\infty) \\ \bx(0) = \bx_0, \,  \bx(\infty) = \bzero}} \! \frac{1}{2} \int\displaylimits_{0}^{\infty} \Vert \by(t) \Vert^2  +  \frac{\Vert \bu(t) \Vert^2}{\eta} {\rm{d}}t,
        \end{equation}
        whereas if $\gamma>1$, the $\cH_\infty$ future energy is defined as
        \begin{equation} \label{eq:HinftyFutureEnergyDef-b}
                \cE_\gamma^{+}(\bx_0)  \coloneqq \!\!\!\! \min_{\substack{\bu \in L_{2}[0,\infty) \\ \bx(0) = \bx_0 , \, \bx(\infty) = \bzero}} \! \frac{1}{2} \int\displaylimits_{0}^{\infty} \Vert \by(t) \Vert^2  +
                \frac{\Vert \bu(t) \Vert^2}{\eta} {\rm{d}}t.
        \end{equation}
\end{definition}

The following theorem states that the energy functions, which are nominally defined by optimization problems, can be computed as the solutions to HJB PDEs.
\begin{theorem}\cite[Thm. 5.2]{Scherpen1996}
        Assume that the HJB equation
        \begin{equation} \label{eq:Hinfty-Past-HJB}
                \begin{split}
                        0 & =  \frac{\partial \cE_\gamma^{-}(\bx)}{\partial \bx} \bf(\bx) + \frac{1}{2}  \frac{\partial \cE_\gamma^{-}(\bx)}{\partial \bx} \bg(\bx) \bg(\bx)^\top \frac{\partial^\top \cE_\gamma^{-}(\bx)}{\partial \bx} \\
                          & \qquad - \frac{\eta}{2}  \bh(\bx)^\top  \bh(\bx)
                \end{split}
        \end{equation}
        has a solution with $\cE_\gamma^{-}(\bzero) = 0$ such that the quantity
        $ - \bf(\bx) - \bg(\bx) \bg(\bx)^\top \partial^\top \cE_\gamma^{-}(\bx)/\partial \bx $
        is asymptotically stable.
        Then this solution is the past energy function $\cE_\gamma^{-}(\bx)$ from~\cref{eq:HinftyPastEnergyDef}.
        Furthermore, assume that the HJB equation
        \begin{equation} \label{eq:Hinfty-Future-HJB}
                \begin{split}
                        0 & =  \frac{\partial \cE_\gamma^{+}(\bx)}{\partial \bx} \bf(\bx)   - \frac{\eta}{2} \frac{\partial \cE_\gamma^{+}(\bx)}{\partial \bx} \bg(\bx) \bg(\bx)^\top \frac{\partial^\top \cE_\gamma^{+}(\bx)}{\partial \bx} \\
                          & \qquad + \frac{1}{2}\bh(\bx)^\top \bh(\bx)
                \end{split}
        \end{equation}
        has a solution with $\cE_\gamma^{+}(\bzero) = 0$ such that the quantity
        $\left.\bf(\bx) - \eta \bg(\bx) \bg(\bx)^\top \partial^\top \cE_\gamma^{+}(\bx)/\partial \bx\right.$
        is asymptotically stable.
        Then this solution is the future energy function $\cE_\gamma^{+}(\bx)$ from~\cref{eq:HinftyFutureEnergyDef-b} for $\gamma>1$ and from~\cref{eq:HinftyFutureEnergyDef-a} for $\gamma <1$.
\end{theorem}

\subsection{Al'brekht's Method To Solve HJB PDEs}\label{sec:albrekht}
In general, solving the HJB PDEs \cref{eq:Hinfty-Past-HJB,eq:Hinfty-Future-HJB} analytically is not feasible.
In this work, in the interest of developing scalable algorithms for nonlinear balancing, we focus on the approximation technique of Al'brekht.
For the special case when the dynamics \cref{eq:FOM-NL} are analytic, then the energy function solutions to the HJB PDEs are also analytic \cite{Albrekht1961,Lukes1969}.
As such, Al'brekht showed that the energy function solutions to \cref{eq:Hinfty-Past-HJB,eq:Hinfty-Future-HJB} can be expanded in a convergent power series
\begin{align*}
         & \cE_\gamma^-(\bx)
        =  {\cE_\gamma^-}^{(2)}(\bx) + {\cE_\gamma^-}^{(3)}(\bx) + \dots + {\cE_\gamma^-}^{(d)}(\bx) + \dots \nonumber \\
         & \cE_\gamma^+(\bx)
        = {\cE_\gamma^+}^{(2)}(\bx) + {\cE_\gamma^+}^{(3)}(\bx) + \dots + {\cE_\gamma^+}^{(d)}(\bx) + \dots
\end{align*}
where $^{(2)}$ denotes a quadratic term, $^{(3)}$ denotes a cubic term, and so on.
According to Al'brekht, upon inserting polynomial expressions into the HJB PDEs \cref{eq:Hinfty-Past-HJB,eq:Hinfty-Future-HJB}, collecting terms of the same polynomial degree gives an algebraic equation for each energy function component ${\cE_\gamma^-}^{(i)}(\bx)$ and ${\cE_\gamma^+}^{(i)}(\bx)$ for $i=2,\dots,d$.
The equations for the quadratic terms ${\cE_\gamma^-}^{(2)}(\bx)$ and ${\cE_\gamma^+}^{(2)}(\bx)$ yield algebraic Riccati equations, and the remaining components ${\cE_\gamma^-}^{(k)}(\bx)$ and ${\cE_\gamma^+}^{(k)}(\bx)$ for $k=3,\dots,d$ solve linear algebraic systems.

While Al'brekht conceptually introduced the idea of a power-series-based approach to solving HJB PDEs, it was presented in an abstract manner without explicit equations to solve and without scalable software.
In this work, we adopt the Kronecker product representation to write the polynomials explicitly so that we may derive the exact form for the algebraic equations for the energy function polynomial representation.

\subsection{Kronecker Product Definitions and Notation}\label{sec:notation}
The Kronecker product of two matrices $\bA \in \real^{p \times q}$ and $\bB \in \real^{s \times t}$ is the $ps \times qt$ block matrix
\begin{align*}
        \bA \otimes \bB \coloneqq \begin{bmatrix} a_{11}\bB & \cdots & a_{1q}\bB \\
                \vdots    & \ddots & \vdots    \\
                a_{p1}\bB & \cdots & a_{pq}\bB
                                  \end{bmatrix},
\end{align*}
where $a_{ij}$ denotes the $(i,j)$th entry of $\bA$.
Repeated Kronecker products are written as
\begin{equation*}
        \kronF{\bx}{k} \coloneqq \underbrace{\bx \otimes \dots \otimes \bx}_{k \ \text{times}}
        \in \real^{n^k}.
\end{equation*}

\noindent For $\bA \in \real^{p \times q}$, the
\textit{$k$-way Lyapunov matrix} is defined as
\small
\begin{equation*}
        \cL_k(\bA) \coloneqq \sum_{i=1}^k\underbrace{\bI_p \otimes \bA \otimes \bI_p \otimes \dots \otimes \bI_p}_{\text{$k$ factors, $\bA$ in the $i$th position}} \in \real^{p^k \times p^{k-1}q}.
\end{equation*}

\normalsize

We also use the $\Vec{[\cdot]}$ operator, which stacks the columns of a matrix into a single column vector, and the \textit{perfect shuffle matrix} $\perm{p}{q} \in \real^{pq \times pq}$ \cite{Henderson1981,VanLoan2000}, defined as the permutation matrix which shuffles $\Vec[\bA]$ to match $\Vec[\bA^\top]$:
\begin{equation*}
        \Vec[\bA^\top]=\perm{q}{p} \Vec[\bA].
\end{equation*}

A concept which arises when dealing with polynomials in Kronecker product form is symmetry of the coefficients (a generalization of symmetry of a matrix), as defined next.
\begin{definition}[Symmetric Coefficients\label{def:sym}] Given a monomial of the form $\bw_d^\top \kronF{\bx}{d}$, the coefficient $\bw_k \in \real^{n^k \times 1}$ is \emph{symmetric} if for all $\ba_i \in \real^n$ it satisfies
        \begin{displaymath}
                \bw_k^\top \left(\ba_1 \otimes \ba_2 \otimes \cdots \otimes \ba_k\right) = \bw_k^\top \left(\ba_{i_1} \otimes \ba_{i_2} \otimes \cdots \otimes \ba_{i_k}\right),
        \end{displaymath}
        where the indices $\{ i_j \}_{j=1}^k$ are any permutation of $\{1, \dots, k \}$.
\end{definition}


\section{Computing \texorpdfstring{$\cH_\infty$}{H-infinity} Energy Function Approximations for Polynomial Drift Systems} \label{sec:NLBT-Poly}
We restrict our consideration to nonlinear dynamical systems with polynomial drift and linear inputs and outputs:
\begin{equation}\label{eq:FOM-Poly}
        \dot{\bx}  = \underbrace{\bA \bx + \sum_{p=2}^\ell \F_p \kronF{\bx}{p}}_{\bf(\bx)} + \underbrace{\bB}_{\bg(\bx)} \bu, \qquad
        \by        = \underbrace{\bC \bx}_{\bh(\bx)},
\end{equation}
where $\bA \in \real^{n\times n}$, $\F_p \in \real^{n\times n^p}$, $\bB \in \real^{n\times m}$, and $\bC \in \real^{p\times n}$.
We seek to compute energy function approximations that solve the HJB PDEs \cref{eq:Hinfty-Past-HJB,eq:Hinfty-Future-HJB} for systems with the particular structure in \cref{eq:FOM-Poly}.
According to Al'brekht's method, the energy function solutions can be represented as polynomials; we truncate the approximations to degree~$d$ and write the energy functions using the Kronecker product.
Defining ${\cE_\gamma^-}^{(i)}(\bx) \coloneqq \frac{1}{2}\bv_i^\top \kronF{\bx}{i}$ and ${\cE_\gamma^+}^{(i)}(\bx) \coloneqq \frac{1}{2}\bw_i^\top \kronF{\bx}{i}$, the energy functions take the explicit forms
\begin{align}
         & \cE_\gamma^-(\bx)
        \approx \frac{1}{2} \sum_{i=2}^d \bv_i^\top \kronF{\bx}{i},
         & \cE_\gamma^+(\bx)
        \approx \frac{1}{2} \sum_{i=2}^d \bw_i^\top \kronF{\bx}{i}, \label{eq:vi-coeffs}
\end{align}
with the coefficients $\bv_i,\bw_i \in \real^{n^i}$ for $i=2, 3, \dots, d$.
In the remainder of this article, the approximations in \cref{eq:vi-coeffs} are treated as equalities.
From here, we derive the algebraic equations for the coefficients $\bv_i$ and $\bw_i$ for $i=2, 3, \dots, d$.
Inserting the polynomial expression \cref{eq:vi-coeffs} for the past energy function $\cE_\gamma^-(\bx)$ into the HJB PDE \cref{eq:Hinfty-Past-HJB}, along with the polynomial dynamics in \cref{eq:FOM-Poly}, the collection of degree~2 terms is
\begin{equation*}
        \begin{split}
                0  = & \frac{1}{2}\left (\bv_2^\top (\bI_n \otimes \bx) + \bv_2^\top (\bx \otimes \bI_n)  \right) \bA \bx                                        \\
                     & + \frac{1}{8}\left (\bv_2^\top (\bI_n \otimes \bx) + \bv_2^\top (\bx \otimes \bI_n)  \right) \bB \times                                   \\
                     & \bB^\top \left (\bv_2^\top (\bI_n \otimes \bx) + \bv_2^\top (\bx \otimes \bI_n)  \right)^\top - \frac{\eta}{2} \bx^\top \bC^\top \bC \bx.
        \end{split}
\end{equation*}
This is a \emph{quadratic} algebraic equation.
Assuming the equation is symmetric\footnote{$\bx^\top \bM\bx$ only implies $\bM = \bzero$ if $\bM$ is symmetric.}, we require equality to hold for all $\bx$; this leads to an algebraic Riccati equation, as shown in \cite{Kramer2022a}.
Al'brekht's method continues by collecting terms of the next degree in a recursive fashion.
The collection of degree~3 and degree~4 terms lead to the \emph{linear} algebraic equations
\begin{align*}
        \cL_{3} \left(\bA + \bB \bB^\top  \bV_2 \right)^\top \mathbf{\tv}_3
         & = - \cL_2(\F_2)^\top \bv_2,
        \\
        \cL_{4} \left(\bA + \bB \bB^\top  \bV_2 \right)^\top \mathbf{\tv}_4
         & =  - \cL_2(\F_3)^\top \bv_2 - \cL_3(\F_2)^\top \bv_3
        \\
         & \hspace{.75cm} - \frac{9}{4} ~\Vec[\bV_3^\top \bB \bB^\top \bV_3].
\end{align*}
Carrying out the process for the remaining coefficients, we arrive at the next two theorems, which give the explicit equations to compute the polynomial coefficients $\bv_i$ and $\bw_i$ for the energy functions~\cref{eq:vi-coeffs} for the polynomial drift dynamics~\cref{eq:FOM-Poly}.
\begin{theorem}[Past energy polynomial coefficients]\label{thm:viPoly}
        Let $\gamma > \gamma_0 \geq 0$ and $\eta = 1-\gamma^{-2}$, where $\gamma_0$ denotes the smallest $\tilde{\gamma}$ such that a stabilizing controller exists for which the $\mathcal{H}_\infty$ norm of the closed-loop system is less than $\tilde{\gamma}$.
        Let the past energy function $\cE_\gamma^{-}(\bx)$, which solves the $\cH_\infty$ HJB PDE \cref{eq:Hinfty-Past-HJB} for the polynomial drift system \cref{eq:FOM-Poly}, be of the form~\cref{eq:vi-coeffs} with the coefficients $\bv_i \in \real^{n^i}$ for $i=2,3,\dots,d$.
        Then $\bv_2 = \Vec\left[\bV_2\right]$, where $\bV_2$ is the symmetric positive semidefinite solution to the $\cH_\infty$ algebraic Riccati equation
        \begin{align}\label{eq:V2-2}
                \bzero & = \bA^\top \bV_2 + \bV_2 \bA- \eta\bC^\top \bC + \bV_2 \bB \bB^\top \bV_2.
        \end{align}
        For $3 \leq k \leq d$, let $\mathbf{\tv}_k \in \real^{n^k}$ solve the linear system
        \small
        \begin{equation}\label{eq:LinSysForVk-Poly}
                \begin{split}
                        \cL_{k} \left(\bA + \bB \bB^\top  \bV_2 \right)^\top \mathbf{\tv}_k
                         & =
                        -\sum_{\mathclap{\substack{i,p\geq 2         \\       i + p = k+1}}}
                        \cL_i(\F_p)^\top \bv_i
                        \\
                         & \hspace{-.75cm}
                        - \frac{1}{4}\sum_{\mathclap{\substack{i,j>2 \\ i+j=k+2}}} ij~\Vec[\bV_i^\top \bB \bB^\top \bV_j]
                \end{split}
        \end{equation}
        \normalsize
        Then the coefficient vector $\bv_k = \Vec\left[\bV_k\right] \in \real^{n^k}$ for $3 \leq k \leq d$ is obtained by symmetrization of $\mathbf{\tv}_k$.
\end{theorem}
\begin{theorem}[Future energy polynomial coefficients]\label{thm:wiPoly}
        Let $\gamma > \gamma_0 \geq 0$ and $\eta = 1- \gamma^{-2}$ as in \cref{thm:viPoly}.
        Let the future energy function $\cE_\gamma^{+}(\bx)$, which solves the $\cH_\infty$ HJB PDE \cref{eq:Hinfty-Future-HJB} for the polynomial drift system \cref{eq:FOM-Poly}, be of the form~\cref{eq:vi-coeffs} with the coefficients $\bw_i \in \real^{n^i}$ for $i=2,3,\dots,d$.
        Then $\bw_2 = \Vec\left[\bW_2\right]$, where $\bW_2$ is the symmetric positive semidefinite solution to the $\cH_\infty$ algebraic Riccati equation
        \begin{align}
                \bzero & = \bA^\top \bW_2 + \bW_2 \bA + \bC^\top \bC - \eta \bW_2 \bB \bB^\top \bW_2.
        \end{align}
        For $3 \leq k \leq d$, let $\mathbf{\tw}_k \in \real^{n^k}$ solve the linear system
        \small
        \begin{equation}\label{eq:LinSysForWk-Poly}
                \begin{split}
                        \cL_{k} \left(\bA - \eta \bB \bB^\top  \bW_2 \right)^\top
                        \mathbf{\tw}_k & =
                        - \sum_{\mathclap{\substack{i,p\geq 2           \\ i + p = k+1}}} \cL_i(\F_p)^\top \bw_i
                        \\ & \hspace{-.75cm}
                        + \frac{\eta}{4}\sum_{\mathclap{\substack{i,j>2 \\ i+j=k+2}}} ij~\Vec[\bW_i^\top \bB \bB^\top \bW_j]
                \end{split}
        \end{equation}
        \normalsize
        Then the coefficient vector $\bw_k = \Vec\left(\bW_k\right) \in \real^{n^k}$ for $3 \leq k \leq d$ is obtained by symmetrization of $\mathbf{\tw}_k$.
\end{theorem}

It is shown in \cite{Kramer2022a} that the matrices $\cL_{k} \left(\bA + \bB \bB^\top  \bV_2 \right)^\top$ and $\cL_{k} \left(\bA - \eta \bB \bB^\top  \bW_2 \right)^\top$ are invertible, and hence the linear systems \cref{eq:LinSysForVk-Poly,eq:LinSysForWk-Poly} are uniquely solvable for the energy function coefficients $\bv_k$ and $\bw_k$.


\section{Implementation Details}\label{sec:complexity}
With a naive implementation, the Kronecker product-based approach scales poorly:
it requires solving linear systems \cref{eq:LinSysForVk-Poly,eq:LinSysForWk-Poly} of dimension $n^k$ for $\bv_k$, which has a cost of $O(n^{3k})$ using a direct solver.
Thus efficiently solving---and even forming---the linear systems \cref{eq:LinSysForVk-Poly,eq:LinSysForWk-Poly} is critical for scalability to systems with large state dimension $n$.
In this section, we outline some of the key details for efficiently implementing the proposed approach, focusing on \cref{eq:LinSysForVk-Poly} since the cost for \cref{eq:LinSysForWk-Poly} is identical.

\subsection*{Forming the right-hand-sides of the linear systems}
To begin, we must consider how to efficiently assemble the right-hand side of \cref{eq:LinSysForVk-Poly}, starting with the first set of terms:
\begin{align}\label{eq:Fterms}
        -\sum_{\mathclap{\substack{i,p\geq 2 \\ i + p = k+1}}} \cL_i(\F_p)^\top \bv_i.
\end{align}
The matrix $\cL_i(\F_p)$ has dimension $(n^i \times n^k)$, whereas the vector $\bv_i$ is $(n^i \times 1)$,
so the cost of evaluating the Lyapunov product $\cL_i(\F_p)^\top \bv_i$ using naive matrix-vector multiplication is $O(n^{k+i})$ using level-2 BLAS operations.
The dominant cost occurs for the term with $i=k-1$, for a total cost of $O(n^{2k-1})$.
Instead, we exploit the structure of the \emph{i-way Lyapunov matrix} to form these terms more efficiently.
Consider a single term from the sum \cref{eq:Fterms}, and expand it according to the definition of the \emph{i-way Lyapunov matrix}
\begin{align*}
        \cL_i(\F_p)^\top \bv_i = (\F_p^\top \otimes \bI_{n^{i-1}}) \bv_i + (\bI_n \otimes \F_p^\top \otimes \bI_{n^{i-2}}) \bv_i + \dots
\end{align*}
All of the quantities on the right-hand side are equivalent under an appropriate permutation/reshaping, so the total cost of forming this is $i$ times the cost of computing the first term.
Using the Kronecker-vec relation \cite{Henderson1981} $(\bD^\top \otimes \bI_p) \Vec[\bA] =  \Vec[\bA \bD]$, we rewrite the first term in the sum as
\begin{align*}
        (\F_p^\top \otimes \bI_{n^{i-1}}) \bv_i & = \Vec[\bV_i^\top \F_p],
\end{align*}
which is now matrix multiplication of $(n^{i-1} \times n)$ and $(n \times n^p)$ matrices, which has a cost of $O(n^{i+p})$ using level-3 BLAS operations.
Since $i+p = k+1$, this is equivalent to $O(n^{k+1})$.
Performing this operation $i$ times for the remaining terms, the total cost of evaluating $\cL_i(\F_p)^\top \bv_i$ this way is $O(i n^{k+1})$.
Again, the dominant cost occurs for the case $i = k - 1$, hence the total cost to form the set of terms \cref{eq:Fterms} is $O(k n^{k+1})$, as opposed to $O(n^{2k-1})$ for a naive implementation.

In the remaining terms in \cref{eq:LinSysForVk-Poly}
\begin{align}
        -\frac{1}{4}\sum_{\mathclap{\substack{i,j>2 \\ i+j=k+2}}} ij~\Vec[\bV_i^\top \bB \bB^\top \bV_j],
\end{align}
products like $\bV_i^\top \bB$ appear repeatedly, so we store them in memory to avoid repeatedly forming them; however, the dominant cost comes from multiplying these stored quantities.
Treating $\bV_i^\top \bB$ as an $(n^{i-1} \times m)$ matrix and $\bB^\top \bV_j$ as an $(m \times n^{j-1})$ matrix, the multiplication $\bV_i^\top \bB \bB^\top \bV_j$ costs $O(m n^{i+j-2})$ using level-3 BLAS operations.
Since $i+j = k+2$, this is equal to $O(m n^k)$.
We form $k-2$ of these terms in the sum, so the overall cost is $O(k m n^k)$.
Assuming $n \gg km$, this is negligible compared to the cost of forming the other terms and solving the equations.

\subsection*{Solving the linear systems}
An efficient solver is required to avoid the $O(n^{3d})$ operations used in a naive direct solve.
As has been shown in other recent works \cite{Breiten2018,Breiten2019,Breiten2019a,Borggaard2020,Borggaard2021,Kramer2022a}, the \emph{k-way Lyapunov matrix} structure of the linear system \cref{eq:LinSysForVk-Poly} can be exploited by modern solvers, e.g., \cite{Grasedyck2004,Chen2019}.
In this work, we use the $k$-way Bartels-Stewart algorithm from \cite{Borggaard2021}, though better performance may be obtained with other solvers.
To summarize the algorithm, a Schur factorization of $\bA + \bB \bB^\top \bV_2$ is used to convert the linear system into an upper-triangular form which requires solving $n^{k-1}$ dense upper triangular linear systems of size $n$.
Each of these $n^{k-1}$ systems can be solved by back substitution with a cost of $O(n^2)$, for a total cost of $O(n^{k+1})$, as opposed to a cost of $O(n^{3k})$ for a naive direct solve of the full system.

\subsection*{Summary of key implementation details}
In summary, the linear system \cref{eq:LinSysForVk-Poly} for the $k$th coefficient $\bv_k$ is formed efficiently by using the Kronecker-vec identity to reshape Kronecker products as matrix-matrix multiplications, which can be computed with $O(kn^{k+1})$ complexity.
Then, using the $k$-way Bartels-Stewart solver from \cite{Borggaard2021}, \cref{eq:LinSysForVk-Poly} can be solved with $O(n^{k+1})$ complexity.
Since the highest-degree coefficient $\bv_d$ is the most expensive to compute, the overall computational complexity is $O(dn^{d+1})$.


\section{Numerical Results for a Nonlinear Heat Equation}\label{sec:fem}
We evaluate the scalability and convergence of the approach proposed in \cref{sec:NLBT-Poly,sec:complexity} on a finite element model of a reaction-diffusion problem.
The results are obtained on a
Linux workstation with an Intel Xeon W-3175X CPU and 256 GB RAM.

\subsubsection*{Model}
Consider the following nonlinear heat equation modeling a reaction-diffusion system
\begin{align*}
        z_t(x,t) & = z_{xx}(x,t) + z_{x}(x,t) + \frac{1}{8} z(x,t)                                                       \\
                 & \hspace{.5cm}  + z(x,t)^3 + \sum_{j=1}^m b_j^m(x) u_j(t)  ,                                           \\
        y_i(t)   & = \int_{\mathrlap{\chi_{[(i-1)\ell/p,i\ell/p]}}} \hspace{.5cm} z(x,t) \text{d} x, \qquad i=1,\dots,p,
\end{align*}
where $z$ represents temperature, $x \in \left[0,\ell\right]$ is the spatial coordinate, and $t$ is time.
The system is subject to Dirichlet boundary conditions
$z(0,t) = z(\ell,t) = 0$ and an initial condition
\begin{align*}
        z(x,0) = 5 \cdot 10^{-5} x(x-\ell)(x-\ell/2).
\end{align*}

Stability of the open-loop system was studied in \cite{Sandstede2005,Galkowski2012,Embree2019}.
Despite the origin being locally asymptotically stable, trajectories that begin near the origin can experience transient growth; this growth feeds the nonlinear terms in the model, making the model unstable despite its linearization being stable.
This feature has been cited as a motivation for nonlinear model reduction methods \cite{Embree2019}.

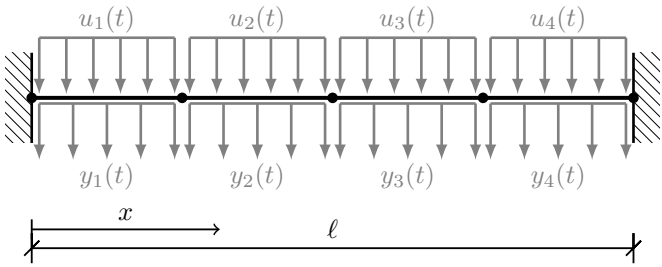
\begin{figure}[htb]
        \centering
        \begin{tikzpicture}
                \fill (0,0)  circle[radius=2pt];
                \fill (2,0)  circle[radius=2pt];
                \fill (4,0)  circle[radius=2pt];
                \fill (6,0)  circle[radius=2pt];
                \fill (8,0)  circle[radius=2pt];

                \point{l}{0}{0};
                \point{r}{8}{0};
                \beam{2}{l}{r};
                \support{3}{l}[-90];
                \support{3}{r}[90];

                \point{ar}{0.1}{-.25};
                \point{br}{2.1}{-.25};
                \point{cr}{4.1}{-.25};
                \point{dr}{6.1}{-.25};

                \point{bl}{1.9}{-.25};
                \point{cl}{3.9}{-.25};
                \point{dl}{5.9}{-.25};
                \point{el}{7.9}{-.25};
                \begin{scope}
                        \tikzstyle{normalLine}=[line width=\normalLineWidth,color=gray]
                        \tikzstyle{smallLine}=[color=white]
                        \tikzstyle{fill}=[color=gray]
                        \lineload{1}{ar}{bl}[.75][.75];
                        \lineload{1}{br}{cl}[.75][.75];
                        \lineload{1}{cr}{dl}[.75][.75];
                        \lineload{1}{dr}{el}[.75][.75];
                \end{scope}

                \point{ar}{0.1}{-1.125};
                \point{br}{2.1}{-1.125};
                \point{cr}{4.1}{-1.125};
                \point{dr}{6.1}{-1.125};

                \point{bl}{1.9}{-1.125};
                \point{cl}{3.9}{-1.125};
                \point{dl}{5.9}{-1.125};
                \point{el}{7.9}{-1.125};

                \begin{scope}
                        \tikzstyle{normalLine}=[line width=\normalLineWidth,color=gray]
                        \tikzstyle{smallLine}=[color=white]
                        \tikzstyle{fill}=[color=gray]
                        \lineload{1}{ar}{bl}[.75][.75][.25];
                        \lineload{1}{br}{cl}[.75][.75][.25];
                        \lineload{1}{cr}{dl}[.75][.75][.25];
                        \lineload{1}{dr}{el}[.75][.75][.25];
                \end{scope}

                \point{x}{2.5}{0};
                \dimensioning{3}{l}{x}{-1.75}[$x$];
                \dimensioning{1}{l}{r}{-2}[$\ell$];

                \node[color=gray] at (1, 1.05) {$u_1(t)$};
                \node[color=gray] at (3, 1.05) {$u_2(t)$};
                \node[color=gray] at (5, 1.05) {$u_3(t)$};
                \node[color=gray] at (7, 1.05) {$u_4(t)$};

                \node[color=gray] at (1, -1.05) {$y_1(t)$};
                \node[color=gray] at (3, -1.05) {$y_2(t)$};
                \node[color=gray] at (5, -1.05) {$y_3(t)$};
                \node[color=gray] at (7, -1.05) {$y_4(t)$};
        \end{tikzpicture}
        \caption{Heat equation on a physical domain of length $\ell$, subdivided here into $m=p=4$ regions for the inputs $u_j(t)$ and outputs $y_i(t)$.}
        \label{fig:heatModel}
\end{figure}

The problem from \cite{Sandstede2005,Galkowski2012,Embree2019} is augmented with control inputs $u_j(t)$ applied over $m$ equally sized subdomains given by the characteristic function $b_j^m(x) = \chi_{[(j-1)\ell/m,j\ell/m]}(x)$ and outputs $y_i(t)$ that are spatial averages of the solution over $p$ equally sized subdomains.
\cref{fig:heatModel} shows a simple diagram of the model.

We pick $\ell = 30$, $m=p=4$, and we discretize the model with $N$ linear finite elements.
This leads to a
finite-dimensional
state-space model
of dimension $n=N-1$
with cubic drift nonlinearity
that can be written as
\begin{align*}
        \dot{\bx}  = \bA \bx + \F_3 \kronF{\bx}{3} + \bB \bu, \qquad
        \by        = \bC \bx.
\end{align*}

\subsubsection*{Results}
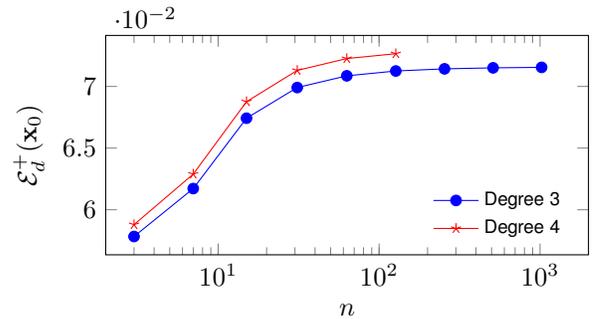
\begin{figure}
        \centering
        \begin{tikzpicture}
                \begin{semilogxaxis}[xlabel=$n$,
                                ylabel=$\mathcal{E}_d^+(\bx_0)$,
                                width=8cm,
                                height=4.5cm,
                                xmin=2,
                                xmax=2000,
                                ylabel near ticks,
                                legend pos=south east,
                                legend style={draw=none,fill=none,font=\sffamily\scriptsize}]

                        \addplot[color=blue,mark=*,mark size=2] table[x={n}, y={E_3^+(x_0)},col sep=&] {example8_convergenceData_d3.dat};
                        \addlegendentry{Degree 3}

                        \addplot[color=red,mark=star,mark size=2] table[x={n}, y={E_4^+(x_0)},col sep=&] {example8_convergenceData_d4.dat};
                        \addlegendentry{Degree 4}
                \end{semilogxaxis}
        \end{tikzpicture}
        \caption{Convergence w.r.t $n$ of the future energy function evaluated at the initial condition $\bx_0$ as the finite element mesh is refined.}
        \label{fig:example8_convergence_n}
\end{figure}

First, we investigate the convergence of the future energy function for $\eta = 0.5$ as the finite element mesh is refined\footnote{For the numerical results, the number of elements $N$ is chosen to be a multiple of 4 so that the boundaries of the control and measurement subdomains occur at nodes in the finite element mesh.}.
Due to spatial limitations, we omit the past energy function results and report that they are qualitatively and quantitatively comparable.
In \cref{tab:example8_convergenceData}, the size of model $n=N-1$ is increased while keeping the degree of the energy function approximation fixed at $d=3$ and $d=4$, respectively.
The energy function values are shown in the second and third columns of \cref{tab:example8_convergenceData} for the initial condition $\bx_0$ corresponding to $z(x,0)$.
In \cref{fig:example8_convergence_n} we plot these energy function values to more clearly show the trend.
As the mesh is refined, the energy function values converge, indicating that the energy function approximations are behaving as expected.

\pgfplotstableread[col sep=&]{example8_convergenceData_d3.dat}{\tableDataA}

\pgfplotstableread[col sep=&]{example8_convergenceData_d4.dat}{\tableDataB}

\pgfplotstablecreatecol[copy column from table={\tableDataB}{[index] 3}] {par1} {\tableDataA}
\pgfplotstablecreatecol[copy column from table={\tableDataB}{[index] 4}] {par2} {\tableDataA}

\makeatletter
\def\pgfmath@error#1#2{\vrule depth 1in}

\def\numorNaN#1{%
        \setbox0\hbox{#1}%
        \ifdim\dp0<1in\box0 \else\makebox[11em]{--}\fi}
\def\bracketorNaN#1{%
        \setbox0\hbox{#1}%
        \ifdim\dp0<1in(\box0 )\else\null\fi}

\begin{table}[htb]
        \centering
        \caption{Future energy function convergence w.r.t. $n$ for $d = 3,4$.}
        \pgfplotstabletypeset[ 
                every head row/.style={before row=\toprule,after row=\midrule},
                every last row/.style={after row=\bottomrule},
                col sep=&,
                multicolumn names,
                columns/{n}/.style={column name=$n$},
                columns/mixed1/.style={string type, column type = {l},column name={$\mathcal{E}_3^+(\bx_0)$ (CPU Sec)}},
                columns/mixed2/.style={string type, column type = {l},column name={$\mathcal{E}_4^+(\bx_0)$ (CPU Sec)}},
                columns={{n},{mixed1},{mixed2}},
                create on use/mixed1/.style={
                                create col/assign/.code={%
                                                \edef\entry{\noexpand\numorNaN{\noexpand\pgfmathprintnumber[sci,sci zerofill, precision=5]{\thisrow{E_3^+(x_0)}}} \noexpand\bracketorNaN{\noexpand\pgfmathprintnumber[sci,sci zerofill, precision=1]{\thisrow{CPU-sec}}}}%
                                                \pgfkeyslet{/pgfplots/table/create col/next content}\entry
                                        }
                        },
                create on use/mixed2/.style={
                                create col/assign/.code={%
                                                \edef\entry{\noexpand\numorNaN{\noexpand\pgfmathprintnumber[sci,sci zerofill, precision=5]{\thisrow{par2}}} \noexpand\bracketorNaN{\noexpand\pgfmathprintnumber[sci,sci zerofill, precision=1]{\thisrow{par1}}}}%
                                                \pgfkeyslet{/pgfplots/table/create col/next content}\entry
                                        }
                        },
        ]{\tableDataA}
        \label{tab:example8_convergenceData}
\end{table}

We confirmed that the energy functions for this initial condition converge with respect to $d$ by computing up to a $d=10$ approximation for $n=7$.
For initial conditions close enough to the origin, the theory is well established that the Taylor approximation converges to the true energy function \cite{Albrekht1961}, so we omit those results since the numerical scaling with respect to $n$ is more of interest for model reduction.

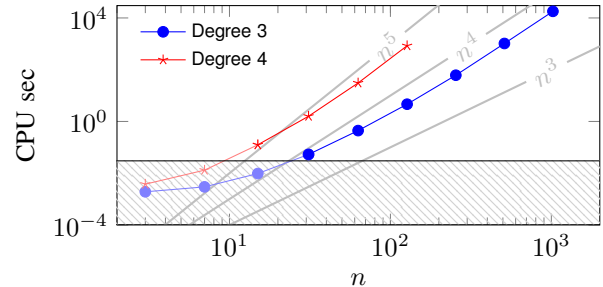
\begin{figure}
        \centering
        \begin{tikzpicture}
                \begin{loglogaxis}[xlabel=$n$,
                                ylabel=CPU sec,
                                width=8cm,
                                height=4.5cm,
                                xmin=2,
                                xmax=2000,
                                ymin=1e-4,
                                ymax=3e4,
                                legend pos=north west,
                                legend style={draw=none,fill=none,font=\sffamily\scriptsize},
                        ]
                        \addplot[domain=1:2e3, samples=10, color=lightgray, thick, forget plot] {10^(-7) * x^3} node [pos=.9, fill=white, sloped] {$n^{3}$};
                        \addplot[domain=1:2e3, samples=10, color=lightgray, thick, forget plot] {10^(-7) * x^4} node [pos=.75, fill=white, sloped] {$n^{4}$};
                        \addplot[domain=1:2e3, samples=10, color=lightgray, thick, forget plot] {10^(-7) * x^5} node [pos=.6, fill=white, sloped] {$n^{5}$};

                        \addplot[color=blue!50,mark=*,mark size=2,skip coords between index={4}{100},forget plot] table[x={n}, y={CPU-sec},col sep=&] {example8_convergenceData_d3.dat};
                        \addplot[color=red!50,mark=star,mark size=2,skip coords between index={3}{100},forget plot] table[x={n}, y={CPU-sec},col sep=&] {example8_convergenceData_d4.dat};

                        \addplot[color=blue,mark=*,mark size=2,skip coords between index={0}{3}] table[x={n}, y={CPU-sec},col sep=&] {example8_convergenceData_d3.dat};
                        \addlegendentry{Degree 3}

                        \addplot[color=red,mark=star,mark size=2,skip coords between index={0}{2}] table[x={n}, y={CPU-sec},col sep=&] {example8_convergenceData_d4.dat};
                        \addlegendentry{Degree 4}

                        \addplot [name path=A,domain=1e-1:2e3,opacity=0] {3e-2};
                        \addplot [name path=B,domain=1e-1:2e3,opacity=0] {1e-4};
                        \addplot [pattern=north west lines, pattern color = gray!35] fill between [of=A and B];
                \end{loglogaxis}
        \end{tikzpicture}
        \caption{Scaling of CPU time w.r.t $n$ for $d=3,4$.
                The scaling closely approaches $O(dn^{d+1})$ as $n$ grows, as predicted by the computational complexity analysis. The grayed out data points are very sensitive to computational variance on the order of tens of milliseconds, so they are not indicative of the true scaling of the algorithms for large models.}
        \label{fig:example8_cpuScaling_n}
\end{figure}

The CPU time required to compute the energy approximations is also shown in \cref{tab:example8_convergenceData}
to demonstrate the scalability of the proposed algorithm
with respect to model dimension $n$.
Since the variance in CPU time can be large for small models (for which the program runs in a fraction of a second), the compute time is averaged over 10 runs for the models of size $n=63$ and smaller.
The computation times are plotted in \cref{fig:example8_cpuScaling_n} to more clearly show the scaling.
As discussed in \cref{sec:complexity}, our algorithm theoretically scales as $O(dn^{d+1})$.
In practice, this is in fact roughly the scaling we see for moderately sized problems, as shown in \cref{fig:example8_cpuScaling_n}.
We observe an apparent departure from the predicted scaling for small $n$, but we attribute this to the minimum overhead time to run the program, which does not scale with the model dimension and appears to be on the order of a few milliseconds.
We also found that the execution time could vary by tens of milliseconds for these small model sizes, so despite averaging over several runs, these data points are highly variable.
We have therefore grayed out the first few data points so that they are still visible but do not obstruct the overall scaling trend as the model dimension increases.

Since the system is odd, the energy functions are even for this example.
Therefore $\bw_3 = \bzero$, and the degree~$3$ energy function approximation is identical to the degree~$2$ LQR solution given by $\bW_2$;
a practitioner seeking to incorporate nonlinear model features in the energy function should therefore compute at least a degree~$4$ approximation.


\section{Conclusion}\label{sec:conclusion}
In this paper, we have presented a novel method for computing $\cH_\infty$ energy functions for nonlinear systems with polynomial drift, linear inputs, and linear outputs.
The approach is based on the Kronecker product polynomial representation, which allows us to obtain explicit formulas for the coefficients of the Taylor-series expansions of the energy functions.
The main advantage over existing approaches is the ability to handle an arbitrary degree of polynomial drift nonlinearity.
By exploiting the tensor structure of the Kronecker product, we have shown that the method is scalable to moderately sized problems, as demonstrated on the nonlinear heat equation example.

The proposed method opens up several directions for future research.
A natural step forward involves extending the approach to nonlinear inputs and outputs so that general polynomial control-affine systems can be studied.
Another direction of major interest is to use the computed energy functions for control and model reduction.
The energy functions can be used directly for polynomial state-feedback controllers, which have shown some promise in the literature over LQR controllers \cite{Borggaard2021,Almubarak2019,Garrard1977}.
The energy functions are also a key step towards obtaining reduced-order models, which can then be used to build observers and controllers for output-feedback problems.
Further accelerating the computations by using more efficient solvers is also of interest.


\bibliography{root}
\bibliographystyle{IEEEtran}

\end{document}